\documentclass[3p]{elsarticle}
\usepackage{amsmath,amsfonts,amssymb}
\usepackage{nicematrix,amsthm}
\usepackage{hyperref}
\usepackage{tikz,mathrsfs}
\biboptions{numbers,sort&compress}
\usetikzlibrary{fit}
\usepackage{textcomp}
\usepackage{float}
\usepackage{placeins}
\usepackage{mathrsfs,mathtools}
\usepackage{algorithm}
\usepackage{algorithmicx,algcompatible}
\usepackage[compatible]{algpseudocode}

\journal{Journal of \LaTeX\ Templates}

\newtheorem{theorem}{Theorem}

\newtheorem{example}{Example}
\newtheorem{lemma}{Lemma}
\newdefinition{remarks}{Remark}
\newdefinition{fact}{Fact}
\newdefinition{definition}{Definition}

\DeclareMathOperator{\rank}{rank}

\DeclareMathOperator{\row}{Row}
\DeclareMathOperator{\blkdiag}{blk-diag}
\DeclareMathOperator{\esssup}{ess~sup}
\DeclareMathOperator{\loc}{loc}

\newcommand{\R}{{\mathbb{R}}}

\newcommand\myae{\stackrel{\mathclap{\normalfont\mbox{a.e.}}}{=}}

\makeatletter
\newenvironment{breakablealgorithm}
{
	\begin{center}
		\refstepcounter{algorithm}
		\hrule height.8pt depth0pt \kern2pt
		\renewcommand{\caption}[2][\relax]{
			{\raggedright\textbf{\fname@algorithm~\thealgorithm} ##2\par}%
			\ifx\relax##1\relax 
			\addcontentsline{loa}{algorithm}{\protect\numberline{\thealgorithm}##2}%
			\else 
			\addcontentsline{loa}{algorithm}{\protect\numberline{\thealgorithm}##1}%
			\fi
			\kern2pt\hrule\kern2pt
		}
	}{
		\kern2pt\hrule\relax
	\end{center}
}
\makeatother

\begin{document}

\begin{frontmatter}

\title{Existence Conditions for  Functional ODE Observer Design of Descriptor Systems Revisited}

\author[a]{Juhi Jaiswal}
\ead{juhi_1821ma03@iitp.ac.in}

\author[b]{Thomas Berger}\ead{thomas.berger@math.upb.de}

\author[a]{Nutan Kumar Tomar\corref{cor1}}
\ead{nktomar@iitp.ac.in}

\address[a]{Department of Mathematics, Indian Institute of Technology Patna, India}

\address[b]{Universit\"at Paderborn, Institut f\"ur Mathematik, Warburger Str.~100, 33098~Paderborn, Germany}

\cortext[cor1]{Corresponding author}

\begin{abstract}
	This paper is devoted to the problem of designing functional observers for linear time-invariant (LTI) descriptor systems. The observers are realized by using state-space systems governed by ordinary differential equations (ODEs). Available existence results for functional ODE observers in the literature are extended by introducing new and milder sufficient conditions. These conditions are purely algebraic and provided directly in terms of the system coefficient matrices. The proposed observer has an order less than or equal to the dimension of the functional vector to be estimated. The observer parameter matrices are obtained by using simple matrix theory, and the design algorithm is illustrated by numerical examples.
\end{abstract}


\end{frontmatter}


	\section{Introduction}\label{intro}
In this paper, we study linear descriptor systems (also known as systems described by differential-algebraic equations (DAEs)) of the form:
\begin{subequations}\label{dls1}
	\begin{eqnarray}
		E\dot{x}(t) &=& Ax(t)+Bu(t), \label{dls1a} \\
		y(t) &=& Cx(t), \label{dls1b} \\
		z(t) &=& Kx(t), \label{dls1c}
	\end{eqnarray}
\end{subequations}
where $E, A \in \mathbb{R}^{m \times n}$, $B \in \mathbb{R}^{m \times l}$, $C \in \mathbb{R}^{p \times n}$, and $K \in \mathbb{R}^{r \times n}$ are known constant matrices. The matrix polynomial $(\lambda E-A)$, in the indeterminate $\lambda$, is called matrix pencil. It is regular, if $m=n$ and $\det (\lambda E-A)$ is not the zero polynomial. In the present paper, we do not assume that  $(\lambda E-A)$ is regular. Note that descriptor systems~\eqref{dls1} occur naturally when dynamical systems are subject to algebraic constraints, see e.g.~\cite{campbell1982singular,dai1989singular,kunkel2006differential,riaza2008differential} and the references therein.

We call $x(t) \in \mathbb{R}^n$ the semistate vector, $u(t) \in \mathbb{R}^l$ the input vector, $y(t) \in \mathbb{R}^p$ the measured output vector, and $z(t) \in \mathbb{R}^r$ the (unknown output) functional vector. The functional vector $z$ contains those variables which cannot be measured and, therefore, observers are required to estimate them. An observer that estimates $z$ without estimating the whole semistate vector~$x$ is said to be a functional (or partial state) observer for \eqref{dls1}. Moreover, any functional observer becomes a full state observer if $K=I_n$, the identity matrix of order $n$.

The contributions made toward functional observers for linear descriptor systems can be broadly classified into two approaches. In the first approach, the observers, known as DAE or generalized observers, consist of a copy or part of the given descriptor system dynamics along with a linear correction term based on the output error \cite{dai1989singular,minamide1989design,berger2019disturbance,JaisBerg23}. On the other hand, the second approach provides ODE observers whose dynamics are described only by ordinary differential equations (ODEs); see \cite{lan2015robust,darouach2012functional,darouach2017functional,jaiswal2021existence,jaiswal2021acc}. An ODE observer is always preferred because it is easy to implement and can be initialized with an arbitrary condition.

There are many applications where we need to estimate only a part or combination of semistates instead of having information about the whole semistate vector. Some particular examples are feedback control and fault or disturbance detection. Besides the applications' viewpoint, theoretically, one of the main reasons for studying functional observers is to find their existence criterion. It is worth pointing out that the functional observers can be designed with considerably weaker assumptions than full state observers~\cite{trinh2011functional}. To this end, investigating characterizations for the existence of functional observers is an active area of research even in the case of standard state space systems, see e.g. \cite{darouach2022functional}. In our recent works \cite{jaiswal2021existence,jaiswal2021acc}, we have established purely algebraic algorithms for functional ODE observer design for~\eqref{dls1} under some assumptions on the system coefficient matrices. These algorithms provide observers of order  $r$, which is the dimension of the functional vector $z$.  Although the observer existence conditions in~\cite{jaiswal2021existence,jaiswal2021acc} are weaker than in the previous works~\cite{lan2015robust,darouach2012functional,darouach2017functional}, they are still not close to being necessary. The purpose of the present paper is to narrow down this gap further and provide milder sufficient conditions for the existence of functional ODE observers than in the above mentioned works.

The novelty of this paper can be summarized as follows:
\begin{enumerate}[1)]
	\item We rigorously introduce the concept of functional ODE observers by considering the same solution framework for \eqref{dls1} as in \cite{jaiswal2021acc,jaiswal2021existence}. Our approach to defining functional ODE observers is based on the behavioral approach as in \cite{berger2019disturbance}. These observers, cf. Definition \ref{def:observer}, cannot lose track once the error between the estimated and true functional vector is zero. It is important to note that the functional ODE observers designed in \cite{darouach2012functional,darouach2017functional,jaiswal2021existence} estimate the vector $z$ asymptotically in the sense of estimators only.
	
	\item A new existence condition for functional ODE observers for \eqref{dls1} is derived. This new condition is milder than those already available in the literature \cite{lan2015robust,darouach2012functional,darouach2017functional,jaiswal2021existence,jaiswal2021acc}.

	\item We propose a new design approach which provides an observer of order less than or equal to $r$. The design algorithm is numerically efficient and mainly consists of three algebraic steps. First, we separate all the zero semistate variables from the dynamics of \eqref{dls1} by orthogonal transformations. In the second step, a relevant part of the dynamics is mixed with the outputs, thanks to structural conditions used for observer design. Finally, the design algorithm resembles the method from~\cite{jaiswal2021existence}.
	
\end{enumerate}

We use the following notation: $0$ and $I$ stand for appropriate dimensional zero and identity matrices, respectively.
In a block partitioned matrix, all missing blocks are zero matrices of appropriate dimensions. Sometimes, for more clarity, the identity matrix of size $n \times n$ is denoted by $I_n$. The set of complex numbers is denoted by $\mathbb{C}$ and $\bar{\mathbb{C}}^+ := \{s \in \mathbb{C}~|~ {\rm Re}\, s\geq 0\}$. The symbols $A^\top$, $A^+$, $\ker A$, and $\row(A)$ denote the transpose, Moore-Penrose inverse (MP-inverse), kernel, and row space of a matrix $A$, respectively. For any square matrix $A$, we write $A > 0$ ($A < 0$) if, and only if, $A$ is positive (negative) definite. A block diagonal matrix having diagonal elements $A_1,\ldots,A_k$ is represented by $\blkdiag\{A_1,\ldots,A_k\}$. 
$\mathscr{L}^1_{\loc}(\mathbb{R}; \mathbb{R}^n)$ is the set of measurable and locally Lebesgue integrable functions $f:\mathbb{R} \to \mathbb{R}^n$,  and $f \myae g$ means $f,~g \in \mathscr{L}^1_{\loc}(\mathbb{R}; \mathbb{R}^n)$ are equal almost everywhere, \emph{i.e.}, $f(t) = g(t)$ for almost all $t \in \mathbb{R}$. Moreover, the notation $``x(t) \rightarrow 0 \text{ as } t \rightarrow \infty"$ means $ ``\lim_{t\to\infty} \esssup \limits_{ [t, \infty)} ||x(t)|| = 0$."

\section{PRELIMINARIES AND PROBLEM STATEMENT}\label{sec:prelim}
The tuple $(x,u,y,z) $ 
is said to be a solution of \eqref{dls1}, if it belongs to the set
\begin{eqnarray*}
	\mathscr{B} := \{(x,u,y,z) \in \mathscr{L}^1_{\loc}(\mathbb{R}; \mathbb{R}^{n+l+p+r})\mid Ex \in \mathcal{AC}_{\loc}(\mathbb{R};  \mathbb{R}^m)  \text{ and } (x,u,y,z)   \text{ satisfies } \eqref{dls1} \text{for almost all } t \in \mathbb{R} \},
\end{eqnarray*}
where $\mathcal{AC}_{\loc}$ represents the set of locally absolutely continuous functions. Our recent study in \cite{jaiswal2021existence} on the existence conditions for functional observers has been established by using the same solution set $\mathscr{B}$; for more details on the solution theory of rectangular systems of type \eqref{dls1}, we refer to~\cite{berger2014,trenn2013solution,berger2017observability}, where $\mathscr{B}$ is called the behavior of \eqref{dls1}. Throughout the paper, we assume that $\mathscr{B}$ is nonempty, which means that there exists an admissible pair for~\eqref{dls1}, consisting of an admissible initial condition and an input function, see also~\cite{hou1999causal,jaiswal2021necessary}.

In this paper, we consider observers of the form
\begin{subequations}\label{obs}
	\begin{eqnarray}
		\dot w(t) & = & Nw(t) + H \begin{bmatrix} u(t) \\ y(t) \end{bmatrix}, \label{obsA} \\
		\hat{z}(t) & = & Rw(t) + M \begin{bmatrix} u(t) \\ y(t) \end{bmatrix}, \label{obsB}
	\end{eqnarray}
\end{subequations}
where $N\in\R^{q\times q}$, $q\le r$, $H\in\R^{q\times (l+p)}$, $R\in\R^{r\times q}$, $M\in\R^{r\times (l+p)}$. An observer of the form~\eqref{obs} is called ODE observer and the non-negative integer $q$ is called its order. We now exploit the behavior $\mathscr{B}$ to define functional ODE observers for \eqref{dls1}.
\begin{definition}\label{def:observer}
	System \eqref{obs} is said to be a functional ODE observer for \eqref{dls1}, if for every  $(x,u,y,z) \in \mathscr{B}$ there exists $w\in \mathcal{AC}_{\loc}(\R;\R^q)$ and $\hat z\in \mathcal{L}^1_{\loc}(\R;\R^r)$ such that $(w,u,y,\hat z)$ satisfy~\eqref{obs} for almost all $t\in\R$, and for all $w,\hat z$ with this property we have
	\begin{enumerate}[(a)]
		\item $ \hat{z}(t) - z(t)\to 0$ for $t\to\infty$,
		\item if $z(0) = \hat{z}(0)$, then $z \myae \hat{z}$.
	\end{enumerate}
\end{definition}

It is important to note that the observers designed in \cite{darouach2012functional,darouach2017functional,jaiswal2021existence} only satisfy condition (a) in  Definition \ref{def:observer} and condition (b) is not required in these works, thus they are only asymptotic estimators. For more details, we refer to Remark~\ref{rem:cond2} and Example \ref{exp2} in the subsequent sections.

The following fundamental results from matrix theory will be required in deriving the main theorem in the next section.

\begin{lemma}\cite{MatsStyn74}\label{prop:blockrank}
	Let $A ,~B$, and $C$ be any matrices of compatible dimensions. If $A$ has full row rank and/or $C$ has full column rank, then $\rank \begin{bmatrix}
	A & B \\ 0 & C
	\end{bmatrix} = \rank{A} + \rank{C}.$
\end{lemma}

\begin{lemma}\cite{jaiswal2021existence}\label{lm:solvability}
	System $XA = B$ has a solution for $X$ if, and only if, $B = BA^+A$, or equivalently, $\rank \begin{bmatrix}
	A \\ B \end{bmatrix} = \rank A$. Moreover,
	$    X = BA^+ - Z(I-AA^+),$
	where $Z$ is an arbitrary matrix of compatible dimension.
\end{lemma}

\begin{lemma}\label{lm:condb}
	Let $A \in \mathbb{R}^{n \times n}$, $C \in \mathbb{R}^{p \times n}$ and set $\mathcal{O}(A,C) = \begin{bmatrix}
	C^\top & (CA)^\top & (CA^2)^\top & \ldots & (CA^{n-1})^\top \end{bmatrix}^\top$. Then the following three statements are equivalent:
	\begin{enumerate}[(i)]
		\item $\forall\, x\in \mathcal{AC}_{\loc}(\mathbb{R};  \mathbb{R}^n)$ with $\dot x(t) = Ax(t)$:\ $Cx(0) = 0 \implies Cx = 0$.
		\item $\ker \mathcal{O}(A,C) = \ker C$.
		\item $\rank \mathcal{O}(A,C) = \rank C$.
	\end{enumerate}
\end{lemma}	

\begin{proof}
	(i) $\Leftrightarrow$ (ii): This is proved in \cite[Lem.~A.1]{berger2019disturbance}.\\
	(ii) $\Leftrightarrow$ (iii):
	Since $\ker C \subseteq \ker \mathcal{O}(A,C)$ it is clear that both subspaces are equal (i.e.,~(ii) holds) if, and only if, they have the same dimension, which is equivalent to~(iii).
\end{proof}

We conclude this section by providing the following decomposition for a given system \eqref{dls1}.
\begin{lemma}\label{lm:Cdecomposition}
	For $E,A\in\R^{n\times m}$ and $B\in\R^{m\times l}$ there exist orthogonal matrices $U \in \mathbb{R}^{m \times m}$ and $V \in \mathbb{R}^{n \times n}$ such that\\
\begin{equation}\label{Cdecomposition}
				UEV = \NiceMatrixOptions
			{nullify-dots,code-for-last-col = \color{black},code-for-last-col=\color{black}}
			\begin{bNiceMatrix}[first-row,last-col]
				& \Ldots[line-style={solid,<->},shorten=0pt]^{k- \text{ block columns}} \\
				E_{11} & E_{1,{k-1}} & \boxtimes & \ldots & \boxtimes &  \\
				0 & E_{2,{k-1}} & \boxtimes & \ldots & \boxtimes &  \\
				& 0 & E_{k-2}  &  \ddots & \vdots & \\
				&  & \ddots & \ddots & \boxtimes & ~\Vdots[line-style={solid,<->},shorten=0pt]^{(k+1)- \text{ block rows}} \\
				&  &  & 0 & E_1 & \\
				& & & & 0 &
			\end{bNiceMatrix}~~, ~	UAV = \NiceMatrixOptions
			{nullify-dots,code-for-last-col = \color{black},code-for-last-col=\color{black}}
			\begin{bNiceMatrix}[first-row,last-col]
				& \Ldots[line-style={solid,<->},shorten=0pt]^{k- \text{ block columns}} \\
				A_{11} & \boxtimes & \boxtimes & \ldots & \boxtimes &  \\
				A_{21} & \boxtimes & \boxtimes & \ldots & \boxtimes &  \\
				& A_{k-1} & \boxtimes & \ldots & \boxtimes & \\
				&  & \ddots & \ddots & \vdots & ~\Vdots[line-style={solid,<->},shorten=0pt]^{(k+1) -\text{ block rows}} \\
				&  & & A_2 & \boxtimes & \\
				& & &  & A_1 & \\
			\end{bNiceMatrix}~~~, ~    UB = \begin{bNiceMatrix}[last-col]
				B_{11} & \\ B_{21} & \\0 & \\ \vdots &  ~\Vdots[line-style={solid,<->},shorten=0pt]^{(k+1) -\text{ block rows}} \\ 0 &
			\end{bNiceMatrix}~~,
		\end{equation}
	where $\boxtimes$ represents the matrix elements of no interest and for each $i=1,~2,~\ldots,~k-1$, where $k \leq n$,
	\begin{enumerate}[(a)]
		\item $A_i$ has full column rank,
		\item $\begin{bmatrix}
		\tilde E_i & \tilde B_i
		\end{bmatrix}$ has full row rank, say $n_i$, where
		$ \tilde E_i = 	\begin{bNiceMatrix}
		E_{11} & E_{1,{k-1}} &  & \boxtimes   \\
		0 & E_{2,{k-1}} &  & \boxtimes  \\
		& \ddots & \ddots &  \\
		& & 0 &    E_i \\
		\end{bNiceMatrix} \text{ and }
		\tilde B_i=\begin{bmatrix} B_{11} \\ B_{21} \\ 0\\ \vdots \\ 0 \end{bmatrix}$, 
		\item $\begin{bmatrix}
		E_{11} & B_{11} \\ 0 & B_{21}
		\end{bmatrix}$ has full row rank, and
		\item $E_{11}$ has full row rank.
	\end{enumerate}
\end{lemma}

The proof of the above lemma is analogous to \cite[Lem.~1]{jaiswal2021necessary} and, hence, omitted. However, for the sake of completeness, a complete method to obtain the matrices $U$ and $V$ is summarized in Algorithm \ref{alg1} below.

\begin{breakablealgorithm}\caption{To compute $U$ and $V$ in Lemma \ref{lm:Cdecomposition}}\label{alg1}
	\begin{enumerate}
		\item Initialize $U_O = I_m$ and $V_O = I_n$.
		\item Find an orthogonal matrix $U_1$ such that $\begin{bmatrix}
		E & B \end{bmatrix}$ is row compressed, \emph{i.e.}, $U_1 \begin{bmatrix}
		E & B \end{bmatrix} = \begin{bmatrix} \hat{E} & \hat{B}_2 \\ 0 & 0 \end{bmatrix}$, where $\begin{bmatrix} \hat{E} & \hat{B}_2 \end{bmatrix}$ has full row rank.
		\item Denote $U_1A = \begin{bmatrix} \check{A}_1 \\ \tilde{A}_1 \end{bmatrix}$, where $\check{A}_1$ has the same number of rows as $\hat{E}$.
		\item Find an orthogonal matrix $V_1$ such that $\tilde{A}_1$ is column compressed on the right side, \emph{i.e.}, $\tilde{A}_1V_1 = \begin{bmatrix} 0 & A_1 \end{bmatrix}$.
		\item Denote $\hat{E}V_1 = \begin{bmatrix} \hat{E}_2 & \boxtimes_1
		\end{bmatrix}$ and $\check{A}_1V_1 = \begin{bmatrix} \hat{A}_2 & \boxtimes
		\end{bmatrix}$, where $\boxtimes_1$ and $\boxtimes$ represent matrices of no interest, having the same number of columns as $A_1$.
		\item Update $U_O := \blkdiag\{U_1,I\}U_O$ and $V_O := V_O\blkdiag\{V_1,I\}$.
		\item Apply the Steps $2$ to $6$ to the triple $(\hat{E}_2,\hat{A}_2,\hat{B}_2)$ $k$- times, until the triple  $(\hat{E}_k,\hat{A}_k,\hat{B}_k)$ has the property that $\begin{bmatrix}
		\hat{E}_k & \hat{B}_k \end{bmatrix}$ has full row rank.
		\item Find an orthogonal row compression matrix $P_o$ such that  $P_o\hat{E}_k = \begin{bmatrix}
		E_{11} \\ 0 \end{bmatrix} $, where $E_{11}$ has full row rank.
		\item $U = \begin{bmatrix}
		P_o & 0 \\ 0 & I \end{bmatrix}U_O$ and $V = V_O$.
	\end{enumerate}
\end{breakablealgorithm}

\section{Existence conditions for observers design}\label{sec:existence}
For a given system \eqref{dls1} and $\lambda\in\mathbb{C}$ set
\allowdisplaybreaks
\begin{eqnarray*}
	& \mathscr{E}=\begin{bNiceMatrix}
		E & 0\\
	\end{bNiceMatrix},~ \mathscr{A}=\begin{bNiceMatrix}
		A & B\end{bNiceMatrix},~ \mathcal{A}=\begin{bNiceMatrix}
		0_{(n-1)n \times m} \\ A \end{bNiceMatrix}, & \\
	&  F = \NiceMatrixOptions
	{nullify-dots,code-for-last-col = \color{black},code-for-last-col=\color{black}}
	\begin{bNiceMatrix}[first-row,last-col]
		& \Ldots[line-style={solid,<->},shorten=0pt]^{n \text{-block columns}} \\
		\mathscr E & \mathscr{A} & & & &  \\
		& \mathscr E & \mathscr{A} & & & \\
		&  &  \ddots & \ddots & & ~\Vdots[line-style={solid,<->}]^{n \text{-block rows}} \\
		& & &   \mathscr{E} & \mathscr{A} & \\
		& & & & \mathscr{E} &
	\end{bNiceMatrix}~~,~ \Psi = \begin{bmatrix}
		F & \mathcal{A} & 0 \\ 0 & E & A \\ 0 & C & 0 \\ 0 & 0 & C  \\ 0 & 0 & K \\ 0 & K & 0 \end{bmatrix},  ~ \Gamma = \begin{bmatrix} F & \mathcal{A} & 0 \\ 0 & E & A \\ 0 & C & 0 \\ 0 & 0 & C  \\ 0 & 0 & K \end{bmatrix}, ~  \Omega(\lambda) = \begin{bmatrix} F & \mathcal{A} & 0 \\ 0 & E & A \\ 0 & C & 0 \\ 0 & 0 & C \\ 0 & K & \lambda K \end{bmatrix}, &
\end{eqnarray*}
and introduce the following two rank conditions:
\begin{equation}\label{H1}
\rank \Psi = \rank \Gamma ,
\end{equation}
and
\begin{equation}\label{H2}
\forall~ \lambda \in \bar{\mathbb{C}}^+:\,\,  \rank \Omega(\lambda) =\rank \Gamma .
\end{equation}
The following theorem is the main result of this article.
\begin{theorem}\label{thm:main}
	Under assumptions \eqref{H1} and \eqref{H2} there exists a functional ODE observer of the form~\eqref{obs} for system~\eqref{dls1}.
\end{theorem}

\begin{proof}
	We divide the proof into several steps:
	
	\textbf{Step $1$:} In view of  Lemma \ref{lm:Cdecomposition}, we first remove all the zero semistates from the dynamics of \eqref{dls1}. By writing \eqref{dls1a} in the form of the decomposition  \eqref{Cdecomposition} with $x=V\begin{bmatrix} x_k^\top & x_{k-1}^\top & \ldots & x_1^\top \end{bmatrix}^\top$ according to the block structure, and solving the obtained system from bottom to top, the fact that $A_i$ has full column rank implies $x_i = 0$, for $1 \leq i \leq k-1$. Therefore, \eqref{dls1a} reduces to
	\begin{eqnarray*}
		\begin{bmatrix}
			E_{11} \\ 0
		\end{bmatrix}\dot x_k(t) = \begin{bmatrix}
			A_{11} \\ A_{21}
		\end{bmatrix}x_k(t) + \begin{bmatrix}
			B_{11} \\ B_{21}
		\end{bmatrix}u(t).
	\end{eqnarray*}
	Moreover, if
   \begin{eqnarray}\label{CV}
		CV = \begin{bmatrix}
			C_k & C_{k-1} & \ldots & C_1
		\end{bmatrix},
	\end{eqnarray}
	and
	\begin{eqnarray}\label{KV}
		KV = \begin{bmatrix}
			K_{11} & K_{k-1} & \ldots & K_1 \end{bmatrix},
	\end{eqnarray}
	system \eqref{dls1} can be rewritten as
	\begin{subequations}\label{dls2}
		\begin{eqnarray}
			E_{11}\dot x_k(t) &=& A_{11}x_k(t) + B_{11}u(t), \label{dls2a}\\
			y_1(t) &=& C_{11}x_k(t), \label{dls2b} \\
			z(t) &=& K_{11}x_k(t), \label{dls2c}
		\end{eqnarray}
	\end{subequations}
	where $y_1(t) = \begin{bmatrix} -B_{21}u(t) \\ y(t) \end{bmatrix}$ and $C_{11} = \begin{bmatrix} A_{21} \\ C_k \end{bmatrix}$. Since the vector $z$ in \eqref{dls2} is the same as in \eqref{dls1}, any functional ODE observer for \eqref{dls2} is also one for the system \eqref{dls1}.
	
	\textbf{Step $2$:} In this step, we separate those variables from the functional vector $z$ that are already known from the new output $y_1$. If $\rank \begin{bmatrix}
	K_{11} \\ C_{11} \end{bmatrix} \neq \rank K_{11} + \rank C_{11}$, then there exist matrices $P$ (a permutation matrix), $S_{11}$, $P_1$ and $P_2$ such that $K_{11} = P\begin{bmatrix}
	S_{11} \\ P_1 S_{11}\\ P_2C_{11} \end{bmatrix}$,  $\row (S_{11}) \cap \row (C_{11}) = \{0\}$, $\row (S_{11}) \subseteq \row (K_{11}) $, and $S_{11}$ has full row rank.
	Thus, the system \eqref{dls2} reduces to
	\begin{subequations}\label{dls3}
		\begin{eqnarray}
			E_{11}\dot x_k(t) &=& A_{11}x_k(t) + B_{11}u(t), \label{dls3a}\\
			y_1(t) &=& C_{11}x_k(t), \label{dls3b} \\
			z(t) &=&  P\begin{bmatrix}
				S_{11} \\ P_1 S_{11} \\ P_2C_{11} \end{bmatrix}x_k(t) = P\begin{bmatrix}
				z_1(t) \\ z_2(t) \\ z_3(t)  \end{bmatrix}, \label{dls3c}
		\end{eqnarray}
	\end{subequations}
	where $z_1 = S_{11}x_k$, $z_2 = P_1 z_1$, and $z_3 = P_2 C_{11}x_k = P_2y_1$. Thus, it is sufficient to design a functional ODE observer for the functional vector $z_1$. Note that, if $$\rank \begin{bmatrix}
	K_{11} \\ C_{11} \end{bmatrix} = \rank K_{11} + \rank C_{11},$$ then $S_{11} = K_{11}$ and $P_2$ is an empty matrix.
	
	\textbf{Step $3$:} In this step, we prove that \eqref{H1} and \eqref{H2} hold for system \eqref{dls1} if, and only if, system~\eqref{dls3} satisfies the following two conditions:
	\begin{equation}\label{trans:reduced:H1}
	\rank \Psi_1 = \rank \Gamma_1,
	\end{equation}
	and,
	\begin{equation}\label{trans:reduced:H2}
	\forall ~\lambda \in \bar{\mathbb{C}}^+:\,\,\rank \Omega_1(\lambda) = \rank \Gamma_1,
	\end{equation}
	where
	$\Gamma_1 = \begin{bmatrix}
	E_{11} & A_{11} \\ C_{11} & 0 \\ 0 & C_{11} \\ 0 & S_{11}	\end{bmatrix},\Psi_1 = \begin{bmatrix}
	E_{11} & A_{11} \\ C_{11} & 0 \\ 0 & C_{11} \\ 0 & S_{11} \\ S_{11} & 0
	\end{bmatrix}$, and $\Omega_1(\lambda)  = \begin{bmatrix}
	E_{11} & A_{11} \\ C_{11} & 0 \\ 0 & C_{11} \\ S_{11} & \lambda S_{11}
	\end{bmatrix}.$ Define, with $U$ and $V$ from Lemma~\ref{lm:Cdecomposition}, $$\mathcal{U} := \begin{bNiceMatrix}
	U &  &  \\
	& \Ddots^{(n+1) \text{ times}} & \\
	& & U \\ & & & I_{2p+r}
	\end{bNiceMatrix}\text{ and } \mathcal{V} := \begin{bNiceMatrix}
	\blkdiag\{V,I_l\} &  &  \\
	& \Ddots^{n \text{ times}} & \\
	& & \blkdiag\{V,I_l\} \\ & & & V \\ &&&& V
	\end{bNiceMatrix}.$$ Clearly,
	$$ \rank \Gamma = \rank (\mathcal{U} \Gamma \mathcal{V}) .$$
	We now write the matrix $\Gamma$ in terms of the matrices $E$, $A$, $B$, and obtain all the $(n+4)$-block rows of the matrix $\mathcal{U} \Gamma \mathcal{V}$.
	To simplify the rank of matrix $\mathcal{U} \Gamma \mathcal{V}$, we repeat the following steps for $i = 1$ to $i = k-1$:
\begin{enumerate}[(i)]
	\item Substitute decomposition \eqref{Cdecomposition} in $i^{th}$-row.
	\item Apply Lemma \ref{prop:blockrank} with respect to the full row rank matrix $\begin{bmatrix}
			\tilde{E}_i & \tilde{B}_i
			\end{bmatrix}$ in $i^{th}$-row (due to statement (b) in Lemma \ref{lm:Cdecomposition}).
	\item Substitute decomposition \eqref{Cdecomposition} in $(i+1)^{th}$-row.
	\item Apply Lemma \ref{prop:blockrank} with respect to the full column rank matrix $\tilde{A}_i$ in $i^{th}$-row, where $\tilde{A}_i = \begin{bmatrix}
			A_{i} & \boxtimes & \ldots & \boxtimes  \\
			& \ddots & \ddots & \vdots \\
			& & A_2 & \boxtimes  \\
			& &  & A_1
			\end{bmatrix}$ (due to statement (a) in Lemma \ref{lm:Cdecomposition}).
\end{enumerate}
Therefore, we obtain
\begin{eqnarray}\label{F1}
	\rank \Gamma = \rho_1 + \rank F_1
\end{eqnarray}
where \begin{eqnarray*}
\rho_1 &=& (n_1 + \rank \tilde{A}_1) + (n_2 + \rank \tilde{A}_2) + \ldots + (n_{k-1} + \rank \tilde{A}_{k-1}) \\
&=& (n_1 + \rank A_1) + (n_2 + \rank A_1 + \rank A_2) + \ldots + (n_{k-1} + \rank A_1 + \rank A_2 + \ldots + \rank A_{k-1}) \\
   && \hspace{11cm} \text{(in view of Lemma \ref{prop:blockrank})} \\
    &=& n_1 + (k-1) \rank A_1 + n_2  + (k-2) \rank A_2 + \cdots + n_{k-2} + 2\rank A_{k-2} + n_{k-1} + \rank A_{k-1}
\end{eqnarray*} and
$$F_1 =\rank \setcounter{MaxMatrixCols}{20}
	\begin{bNiceMatrix}[last-col]
		\begin{matrix}
		0 & {A_{k-1}} & { \boxtimes} & {\dots }& { \boxtimes} \\
		\vdots & & {\ddots}  & {\ddots} & {\vdots}  \\
		0 & & & {A_2} & { \boxtimes} \\
		0 &  &  &  & {A_1 }
		\end{matrix} & \begin{matrix}
		0 \\ \vdots \\ 0 \\ 0
		\end{matrix} & & & & & & \\  \hline
		\begin{matrix}
		UEV
		\end{matrix} & \begin{matrix}
		0
		\end{matrix} & \begin{matrix}
		UAV
		\end{matrix} & \begin{matrix}
		BV
		\end{matrix} & & & & & \\
		& & & \ddots & \ddots & & & & &
		~ \Vdots[line-style={solid,<->}]^{(n +5 -k)- \text{ block rows}}
		\\
		& & & \begin{matrix}
		UEV
		\end{matrix} & \begin{matrix}
		0
		\end{matrix} & \begin{matrix}
		UAV
		\end{matrix} & \begin{matrix}
		BV
		\end{matrix} & & \\ \hline
		& & & & & \begin{matrix}
		UEV
		\end{matrix} & \begin{matrix}
		0
		\end{matrix} & \begin{matrix}
		UAV
		\end{matrix} & \\
		& & & & &  & & \begin{matrix}
		UEV
		\end{matrix} & \begin{matrix}
		UAV
		\end{matrix} \\
		& & & & &  & & \begin{matrix}
		CV
		\end{matrix} & \begin{matrix}
		0
		\end{matrix} \\
		& & & & &  & & & \begin{matrix}
		CV
		\end{matrix} \\
		& & & & &  & & & \begin{matrix}
		KV
		\end{matrix} &
		\end{bNiceMatrix}~~.$$
	Now, to simplify the rank of matrix $F_1$, we perform the following steps repeatedly for $i = k$ to $i = n$:
\begin{enumerate}[(i)]
	\item Substitute decomposition \eqref{Cdecomposition} in $(i+1)^{th}$-row.
	\item Apply Lemma \ref{prop:blockrank} with respect to the full column rank matrix $\tilde{A}_{k-1}$ in $i^{th}$-row, where $\tilde{A}_{k-1} = \begin{bmatrix}
		A_{k-1} & \boxtimes & \ldots & \boxtimes  \\
		& \ddots & \ddots & \vdots \\
		& & A_2 & \boxtimes  \\
		& &  & A_1
		\end{bmatrix}$ (due to statement (a) in Lemma \ref{lm:Cdecomposition}).
\item Apply Lemma \ref{prop:blockrank} with respect to the full row rank matrix
$\begin{bmatrix}
E_{11} & B_{11} \\  & B_{21}
\end{bmatrix}$ in $(i+1)^{th}$-row (due to statement (c) in Lemma \ref{lm:Cdecomposition}).		
\end{enumerate}
Therefore, we obtain
\begin{eqnarray}\label{FO}
		\rank \Gamma = \rho_2 + \rank F_2
\end{eqnarray}
where \begin{eqnarray*}
\rho_2 &=&  \rho_1 + (n-k) \rank \tilde{A}_{k-1} + (n-k)\rank\begin{bmatrix}	E_{11} &  B_{11} \\ & B_{21} \end{bmatrix} \\
&=& \rho_1 + (n-k)(\rank A_1 + \rank A_2 + \ldots + \rank A_{k-1} ) + (n-k)\rank\begin{bmatrix}	E_{11} &  B_{11} \\ & B_{21} \end{bmatrix} \\
&& \hspace{11cm} \text{(in view of Lemma \ref{prop:blockrank})} \\
&=& n_1 + (n-1)\rank A_1 + n_2+(n-2)\rank A_2 + ~\ldots~	+n_{k-1}+(n-(k-1)) \rank A_{k-1} \\ && +~(n-k)\rank\begin{bmatrix}	E_{11} &  B_{11} \\ & B_{21} \end{bmatrix}
	\end{eqnarray*} and
	$$F_2 =  \setcounter{MaxMatrixCols}{15}
	\begin{bNiceMatrix}[vlines,rules/width=0.2pt]
	\begin{matrix}
	{E_{11} } \\0 \\ 0 \\ \vdots \\ 0 \\ 0
	\end{matrix}&   \begin{matrix}
	A_{11} & \boxtimes & \boxtimes & \ldots & \boxtimes   \\
	A_{21} & \boxtimes & \boxtimes & \ldots & \boxtimes   \\
	& {A_{k-1}} & {\boxtimes} & {\ldots} & {\boxtimes}  \\
	&  & {\ddots} & {\ddots} & {\vdots}  \\
	&  & & {A_2} & {\boxtimes} \\
	& & &  & {A_1 }
	\end{matrix} &  \\ \hline
	&  \begin{matrix}
	UEV
	\end{matrix}&\begin{matrix}
	UAV
	\end{matrix} \\
	&   \begin{matrix}
	CV
	\end{matrix} & 0 \\
	& &  \begin{matrix} CV \end{matrix} \\
	& &  \begin{matrix}
	KV
	\end{matrix}
	\end{bNiceMatrix}.$$
	We now simplify the rank of $F_2$ in \eqref{FO} by performing the following four steps.
\begin{enumerate}[(i)]
\item Apply Lemma \ref{prop:blockrank} with respect to the full row rank matrix $E_{11}$.
\item Substitute decompositions \eqref{Cdecomposition}, \eqref{CV}, and \eqref{KV} in  the last four block rows in $F_2$.
\item Apply Lemma \ref{prop:blockrank} with respect to the full column rank matrices $\tilde{A}_{k-1}$ in $n^{th}$- and $(n+1)^{th}$-block rows 
\item Substitute $C_{11} = \begin{bmatrix}
		A_{21} \\ C_{k}
		\end{bmatrix}$ and $K_{11} = P\begin{bmatrix}
		S_{11} \\ P_1 S_{11}\\ P_2C_{11}
		\end{bmatrix}$.
	\end{enumerate}
Thus, we obtain
\begin{eqnarray}\label{R1}
	\rank \Gamma =  \rho + \rank \Gamma_1 ,
	\end{eqnarray}
	where
	$\rho = \rho_2 + \rank E_{11} + 2( \rank A_1 + \ldots + \rank A_{k-1} ).$
	
	In a similar manner, it is straightforward to obtain that
	\begin{eqnarray}
		\rank \Psi &=&  \rho + \rank \Psi_1 \label{R2} \\
		\text{and } ~\rank \Omega(\lambda)  &=&  \rho + \rank \Omega_1(\lambda) . \label{R3}
	\end{eqnarray}
	Thus, it follows from \eqref{R1}, \eqref{R2}, and \eqref{R3} that  \eqref{H1} is equivalent to \eqref{trans:reduced:H1} and \eqref{H2} is equivalent to \eqref{trans:reduced:H2}.
	
	\textbf{Step $4$:}  We show that there are matrices $N,~T,~L$, and $\bar{M}$ of appropriate dimensions such that the following system~\eqref{obsv} is a functional ODE observer for system~\eqref{dls3} and, therefore, due to steps $1-2$, is a functional ODE observer for system~\eqref{dls1} as well.
	\begin{subequations}\label{obsv}
		\begin{eqnarray}
			\dot{w}(t) &=& Nw(t) + TB_{11}u(t) + Ly_1(t), \label{obsva} \\
			\hat{z}(t) &=& P(Rw(t) + M y_1(t)), \label{obsvb}
		\end{eqnarray}
	\end{subequations}
	where $q = \rank{S_{11}}$, $R = \begin{bmatrix} I_q \\ P_1\\ 0 \end{bmatrix}$, and $M = \begin{bmatrix} \bar M \\ P_1 \bar M\\ P_2 \end{bmatrix}$. Set $e = \hat{z} - z$ and $e_1 = w - TE_{11}x_k$, then by~\eqref{dls3} and~\eqref{obsv} we have
	\begin{eqnarray}
		\dot{e}_1(t) &=& \dot{w}(t) - TE_{11}\dot{x}_k(t) \nonumber\\
		&=& Ne_1(t) + (NTE_{11} + LC_{11} - TA_{11})x_k(t), \label{errordya}
	\end{eqnarray}
	and
	\begin{eqnarray}
		e(t) &=& P(Rw(t) +My_1(t)) - K_{11}x_k(t) \nonumber\\
		&=& P\left(Rw(t) + \begin{bmatrix} \bar M \\ P_1\bar M \\ P_2 \end{bmatrix}y_1(t) - \begin{bmatrix} S_{11}x_k(t)\\ P_1 S_{11} x_k(t)\\ P_2y_1(t) \end{bmatrix}\right) \nonumber\\
		&=& PR(w(t) + \bar M y_1(t) - S_{11}x_k(t)) \nonumber \\
		&=& PRe_1(t) + PR(TE_{11} + \bar M C_{11} - S_{11})x_k(t). \label{errordyb}
	\end{eqnarray}
	Thus $e(t) \rightarrow 0$ as $t \rightarrow \infty$, if
	\begin{subequations}\label{solv1:cond}
		\begin{eqnarray}
			NTE_{11} + LC_{11} &=& TA_{11}, \label{solv1:conda} \\
			TE_{11} + \bar{M}C_{11} &=& S_{11},  \label{solv1:condb} \\
			\text{ and $N$ is stable.}&& \label{solv1:condc}
		\end{eqnarray}
	\end{subequations}
	Now, by using \eqref{solv1:condb} in \eqref{solv1:conda}, we obtain that \eqref{solv1:cond} is equivalent to
	\begin{subequations}\label{solv:cond}
		\begin{eqnarray}
			TA_{11} + QC_{11} - NS_{11} &=& 0, \label{solv:conda}\\
			TE_{11} + \bar{M}C_{11} &=& S_{11}, \label{solv:condb} \\
			\text{ and $N$ is stable,}&& \label{solv:condc}
		\end{eqnarray}
	\end{subequations}
	where $Q = N\bar M - L$.
	Thus, we have shown that, if there are matrices $N,~T,~Q$, and $\bar{M}$ which satisfy \eqref{solv:cond}, then the observer \eqref{obsv} meets the property (a) in Definition \ref{def:observer}. To prove that \eqref{obsv} also satisfies property (b) in Definition \ref{def:observer}, we first see that, under \eqref{solv:cond}, the error dynamics \eqref{errordya}--\eqref{errordyb} reduces to
	\begin{equation*}
	\dot{e}_1(t) = Ne_1(t) \text{ and } e(t) = PRe_1(t).
	\end{equation*}
	Then, it follows from Lemma \ref{lm:condb} that
	\begin{equation*}
	\forall\, e_1\in \mathcal{AC}_{\loc}(\mathbb{R};  \mathbb{R}^q) \text{ with } \dot{e}_1(t) = Ne_1(t) \text{ for all $t\in\R$}:~ e(0) = 0 \implies e = 0
	\end{equation*}
	if, and only if,
	\begin{eqnarray}\label{rank1}
		\rank(PR) = \rank \mathcal{O}(N,PR).
	\end{eqnarray}
	Since $P$ is a permutation matrix, the rank identity \eqref{rank1} is equivalent to	
	\begin{eqnarray}\label{rank}
		\rank{R} = \rank \mathcal{O}(N,R) .
	\end{eqnarray}
	Clearly, the condition \eqref{rank} holds trivially, because $R$ has full column rank. Hence, under \eqref{solv:cond} the system \eqref{obsv} satisfies property (b) in Definition \ref{def:observer}.
	
	\textbf{Step $5$:} By step $4$ we have that, if the matrix equation \eqref{solv:cond} is solvable, then \eqref{obsv} is a functional ODE observer for \eqref{dls1}. In this step, we show that conditions \eqref{H1} and \eqref{H2} ensure the solvability of \eqref{solv:cond}. By steps $1-3$, it is sufficient to show that, if system \eqref{dls3} satisfies conditions \eqref{trans:reduced:H1} and \eqref{trans:reduced:H2}, then \eqref{solv:cond} has a solution $N,~T,~Q$, and $\bar{M}$. Clearly, Eqs. \eqref{solv:conda} and \eqref{solv:condb} can be rewritten as
\begin{eqnarray}\label{eq:cond:comb}
	\begin{bmatrix}
		T & \bar{M} & Q & -N
	\end{bmatrix} \Gamma_1 = \begin{bmatrix}
		S_{11} & 0 \end{bmatrix}
	\end{eqnarray}
	and by Lemma \ref{lm:solvability} this equation has a solution if, and only if, \eqref{trans:reduced:H1} holds. Moreover, Lemma \ref{lm:solvability} gives that
	\begin{eqnarray}\label{eq:soln}
		\begin{bmatrix}
			T & \bar{M} & Q & -N
		\end{bmatrix} =	 \begin{bmatrix}
			S_{11} & 0 \end{bmatrix} \Gamma_1^+ - Z (I-\Gamma_1 \Gamma_1^+),\quad
	\end{eqnarray}
	where $Z$ is an arbitrary matrix of appropriate dimension.
	Now, by the same argument as in Step~5 of the proof of Theorem~1 in~\cite{jaiswal2021existence} it follows that under condition~\eqref{trans:reduced:H1} there exists a matrix $Z$ of compatible dimension such that $N$ obtained in~\eqref{eq:soln} is stable if, and only if,~\eqref{trans:reduced:H2} holds.
\end{proof}

The following remarks are warranted in view of the proof of Theorem \ref{thm:main}.

\begin{remarks}\label{rem:static}
	If  $\rank \begin{bmatrix}
	K_{11} \\ C_{11} \end{bmatrix} = \rank C_{11}$, then Lemma~\ref{lm:solvability} implies that there exists a matrix $P_2$ such that $K_{11} = P_2 C_{11}$ and hence step~$2$ of Theorem \ref{thm:main} reveals that there exists a static observer, \textit{i.e.,} an observer of order zero, with:
	\begin{eqnarray*}
		\hat{z}(t) = P_2y_1(t) = P_2\begin{bmatrix}
			-B_{21}u(t) \\ y(t)
		\end{bmatrix}
	\end{eqnarray*}
	and in this case $\hat{z}~ \myae ~z.$
\end{remarks}

\begin{remarks}
	It is well known that for state-space systems $\dot x(t) = Ax(t),\ y(t) = Cx(t)$ there exists an observer if, and only if, the pair $(A,C)$ is detectable. Thus, one may wonder to which extent the assumptions~\eqref{H1} and~\eqref{H2} are related to a detectability condition. As it turns out, in the design of functional observers for linear descriptor systems, the correct notion to consider is \textit{partial detectability}, which was introduced and characterized in our recent work \cite{JaisBerg23}. Roughly speaking, this concept means that the functional vector $z$ can be asymptotically obtained from the knowledge of $u$ and $y$. Moreover, it is a direct consequence of Theorem \ref{thm:main} and \cite[Thm. $2$]{JaisBerg23} that a system~\eqref{dls1} which satisfies \eqref{H1} and \eqref{H2} is also partially detectable.
\end{remarks}

\begin{remarks}\label{rem:checkingcriteria}
	It is not always recommended to verify assumption \eqref{H2} initially as it requires the determination of the finite eigenvalues of the pair $(\bar{E},\bar{A})$, where $\bar{E} = \begin{bmatrix}
	E \\ 0 \end{bmatrix}$ and $\bar{A} = \begin{bmatrix}
	A \\ C \end{bmatrix}$. It is a simple observation from Step $5$ in Theorem \ref{thm:main} that under assumption \eqref{H1}, the condition \eqref{H2} holds if, and only if, the matrix pair $(N_1,~N_2)$ is detectable, where
	\begin{equation}\label{eq:N1N2}
	    N_1 =   \begin{bmatrix} S_{11} & 0 \end{bmatrix} \Gamma_1^+ \begin{bmatrix}
	    0 \\ 0 \\ 0 \\ -I \end{bmatrix},\text{ and } N_2 = (I-\Gamma_1 \Gamma_1^+)\begin{bmatrix}
	    0 \\ 0 \\ 0 \\ -I \end{bmatrix}.
	\end{equation}
Moreover, this matrix pair $(N_1,~ N_2)$ is detectable if, and only if, there exists a matrix $Z$ in \eqref{eq:soln} such that $N = N_1 - ZN_2$ is stable, cf. \cite[Thm. $3.4$]{zhou1998essentials}. Furthermore, the celebrated Lyapunov stability theory \cite{boyd1994linear} yields that $N = N_1 - ZN_2$ is stable  if, and only if, there exists a symmetric matrix $\bar{P} > 0$ such that
\begin{equation}\label{Stable1}
	(N_1 - ZN_2)^\top \bar{P} + \bar{P} (N_1 - ZN_2) < 0.
\end{equation}
If we assume $\bar{Z} = \bar{P}Z$, then the linear matrix inequality (LMI) \eqref{Stable1} is equivalent to
\begin{equation}\label{Stable}
	N^\top_1\bar{P} + \bar{P}N_1 - N^\top_2{\bar{Z}}^\top - \bar{Z}N_2 < 0.
	\end{equation}
Thus, there exists $Z$ in \eqref{eq:soln} such that $N$ is stable if, and only if, the LMI \eqref{Stable} is solvable for $\bar{P}$ and $\bar{Z}$ with $\bar{P} = \bar{P}^\top >0$, and if a solution exists, then $Z = \bar{P}^{-1}\bar{Z}$. Notably, the LMI \eqref{Stable} can be solved in a numerically efficient way by using convex optimization based solvers, for instance, `feasp' in the MATLAB LMI toolbox.

\end{remarks}

\begin{remarks}\label{rem:cond2}
	We emphasize that Algorithm~$2$ in~\cite{jaiswal2021existence} does not guarantee the existence of matrices $R$ and $N$ such that the rank identity~\eqref{rank} is satisfied, cf. Example~\ref{exp2}. To the best of our knowledge, the present work is the first one presenting a unified framework for designing functional ODE observers satisfying the observer matching condition~(b) in Definition~\ref{def:observer}. Clearly, the rank identity \eqref{rank} ensures this condition for the proposed observer \eqref{obsv}.
\end{remarks}	

\begin{remarks}\label{ex:sufficiency}
	Theorem~\ref{thm:main} shows that the conditions~\eqref{H1} and~\eqref{H2} are  sufficient for the existence of functional ODE observers. However, these conditions are not necessary. For example, consider system \eqref{dls1} with the coefficient matrices
	\allowdisplaybreaks
	\begin{eqnarray*}
		& E = \begin{bmatrix}
			1 & 0 & 0 & 0 \\ 0 & 1 & 0 & 0 \\ 0 & 0 & 1 & 0 \\ 0 & 0 & 0 & 0
		\end{bmatrix},~A = \begin{bmatrix}
			-3 & 1 & 0 & 0 \\
			2 & -1 & 0 & 0 \\ 0  & 0 & 1 & 0 \\ 0 & 0  & 0 & 1
		\end{bmatrix},~B = \begin{bNiceMatrix}
			0 \\ 1 \\ 0 \\ 1
		\end{bNiceMatrix}, ~ C = \begin{bmatrix}
			0 & 0 & 1 & 0
		\end{bmatrix}, \text{ and }
		K = \begin{bmatrix}
			1 & 1.366 & 0 & 0
		\end{bmatrix} ,&
	\end{eqnarray*}
	then condition \eqref{H1} does not hold. However, this system can be written as
	\begin{subequations}\label{exp:e1}
		\begin{eqnarray}
			\begin{bmatrix}
				\dot{x}_1(t) \\ \dot{x}_2(t)
			\end{bmatrix}   &=& \begin{bmatrix}
				-3 & 1 \\ 2 & -1
			\end{bmatrix} \begin{bmatrix}
				x_1(t) \\ x_2(t)
			\end{bmatrix} + \begin{bmatrix}
				0 \\ 1
			\end{bmatrix}u(t),\label{exp:e1a} \\
			0 &=& x_4(t) + u(t), \label{exp:e1b}\\
			\dot{x}_3(t) &=& x_3(t), \label{exp:e1c}\\
			y(t) &=& x_3(t), \label{exp:e1d}\\
			z(t) &=& x_1(t) + 1.366 x_2(t), \label{exp:e1e}
		\end{eqnarray}
	\end{subequations}
	and one can easily verify that the following system is a functional ODE observer for~\eqref{exp:e1} in the sense of Definition~\ref{def:observer}:
	\begin{eqnarray*}
		\dot{w}(t) &=& \begin{bmatrix}
			-3 & 1 \\ 2 & -1
		\end{bmatrix}w(t) + \begin{bmatrix}
			0 \\ 1
		\end{bmatrix}u(t), \\
		\hat{z}(t) &=& \begin{bmatrix}
			1 & 1.366
		\end{bmatrix} w(t) .
	\end{eqnarray*}
\end{remarks}

\section{Numerical illustration}\label{sec:numerical}

\begin{example}\label{exp2}
	Consider \eqref{dls1} with the coefficient matrices:
	\begin{eqnarray*}
		& E = \begin{bmatrix} 1 & 0 & 0 \\ 0 & 0 & 0 \\ 0 & 1 & 0 \end{bmatrix},\, A = \begin{bmatrix} -1 & 0 & 0 \\ 0 & 1 & 0 \\ 0 & 0 & 1 \end{bmatrix},\, B = \begin{bmatrix} 0 & 1 \\ 0 & 0 \\ 1 & 0 \end{bmatrix}, ~ C = \begin{bmatrix} 1 & 0 & 0 \end{bmatrix}, \text{ and }
		K = \begin{bmatrix} -1 & 1 & 0 \\ 0 & 1 & 0 \end{bmatrix}, &
	\end{eqnarray*}	
	for which conditions \eqref{H1} and \eqref{H2} are satisfied. Hence, we can design a functional ODE observer for system \eqref{dls1} via our proposed algorithm (in the proof of Theorem~\ref{thm:main}).	After removing the zero semistates, we obtain the reduced system~\eqref{dls2} with coefficient matrices:
	\begin{eqnarray*}
		& E_{11} = \begin{bmatrix}
			-1 & 0  \end{bmatrix},~
		A_{11} = \begin{bmatrix}
			1 & 0  \end{bmatrix}, ~B_{11} = \begin{bmatrix}
			0 & 1 \end{bmatrix},  ~ B_{21} = \begin{bmatrix}
				-1 & 0
			\end{bmatrix},~ C_{11} = \begin{bmatrix} 0 & -1 \\ 1 & 0 \end{bmatrix}, \text{ and }
		K_{11} = \begin{bmatrix} -1 & 0 \\ 0 & 0 \end{bmatrix}.&
	\end{eqnarray*}
	Here, $\rank \begin{bmatrix}
	K_{11} \\ C_{11} \end{bmatrix} = 2 = \rank C_{11}$, hence by Remark \ref{rem:static} we have that there exists a static functional observer:
	\begin{equation}\label{exp1:obsv}
	\hat{z}(t) = \begin{bmatrix}
	0 & -1 \\ 0 & 0
	\end{bmatrix}\begin{bmatrix}
	\begin{bmatrix}
	1 & 0
	\end{bmatrix}u(t) \\y(t)
	\end{bmatrix} = \begin{bmatrix}
	-y(t) \\ 0
	\end{bmatrix}.
	\end{equation}
	
	Since $x_2 = 0$, it is clear that $\hat{z} = z$. For comparison, we have designed an observer for the given system by using the design algorithm from~\cite{jaiswal2021existence} and plotted the error vector $e(t) = \hat{z}(t) - z(t)$ in Figure \ref{fig:fig2} with initial conditions $x_k(0) = \begin{bmatrix}
		1 \\ 2
		\end{bmatrix}$, $w(0) = \begin{bmatrix}
		3 \\ 4
		\end{bmatrix}$ and input $u(t) = \begin{bmatrix}
		\sin(t) \\ \exp(t)
		\end{bmatrix}$. Clearly, $e(2.43) = 0$ but $e(t) \neq 0$, for all $t > 2.43$. Hence, in general, the observer presented in~\cite{jaiswal2021existence} does not satisfy condition~(b) in Definition~\ref{def:observer} and works as an estimator only.
	
	\begin{figure}[h]
		\centering
		\begin{center}
			\includegraphics[scale=0.60]{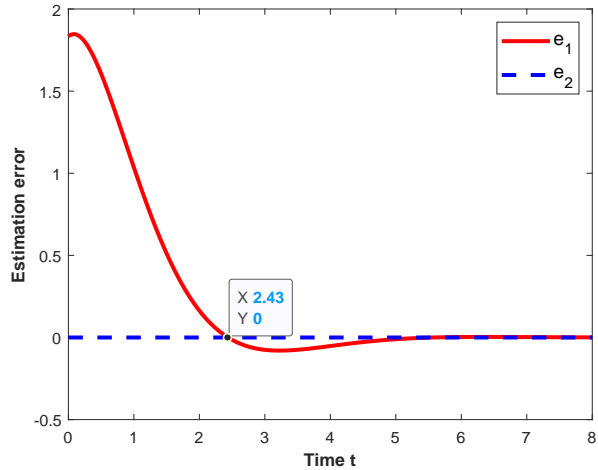}
			\caption{Time response of estimation error obtained from \cite[Algorithm~$2$]{jaiswal2021existence} for Example~\ref{exp2}.}
			\label{fig:fig2}
		\end{center}
	\end{figure}
\end{example}

\begin{example}\label{exp3}
	Consider \eqref{dls1} with the coefficient matrices:
	\begin{eqnarray*}
		& E = \begin{bmatrix} 1 & 0 & 0 & 0 \\ 0 & 0 & 0 & 1 \\ 0 & 1 & 0 & 0 \\ 0 & 0 & 0 & 0 \end{bmatrix}, A = \begin{bmatrix} -1 & 0 & 0 & 0 \\ 0 & 1 & 0 & 0 \\ 0 & 0 & -1 & 0 \\ 0 & 0 & 0 & 1 \end{bmatrix}, B = \begin{bmatrix} 0 \\ 0 \\ 1 \\ 0 \end{bmatrix},~ C = \begin{bmatrix} 0 & 1 & 0 & 0 \end{bmatrix}, \text{ and }
		K = \begin{bmatrix} 1 & 0 & 1 & 1 \end{bmatrix}. &
	\end{eqnarray*}	
	
	This system does not satisfy the assumptions used in any of the research articles \cite{lan2015robust,darouach2012functional,darouach2017functional,jaiswal2021acc,jaiswal2021existence}. However, the system coefficient matrices satisfy conditions \eqref{H1} and \eqref{H2}.
	After removing the zero semistates, the following coefficient matrices are obtained for system \eqref{dls2}:
	\begin{eqnarray*}
		& E_{11} = \begin{bmatrix} 0 & 1 \end{bmatrix},~
		A_{11} = \begin{bmatrix} 0 & -1 \end{bmatrix},~
		B_{11} = \begin{bmatrix} 0 \end{bmatrix},~
		B_{21} = \begin{bmatrix} 1 \end{bmatrix}, ~ C_{11} = \begin{bmatrix} -1 & 0 \\ 0 & 0 \end{bmatrix}, \text{ and }
		K_{11} = \begin{bmatrix} 1 & 1 \end{bmatrix}.&
	\end{eqnarray*}
	Here, $\rank \begin{bmatrix}
	K_{11}\\ C_{11} \end{bmatrix} =  2 = \rank K_{11} + \rank C_{11}$, thus the proposed algorithm provides the following functional ODE observer:
	\begin{eqnarray*}
		\dot{w}(t) &=& -w(t) \\
		\hat{z}(t) &=& w(t) + \begin{bmatrix}
			-1 & 0 \end{bmatrix}y_1(t) = w(t) + u(t).
	\end{eqnarray*}
	In Figure \ref{fig:fig3}, the time responses of $z(t)$ and $\hat{z}(t)$ are plotted by taking $x_k(0) = \begin{bmatrix}
	1 \\ 2 \end{bmatrix},~w(0) = 3$, and $u(t) = \sin(t) $. It is clear from Figure \ref{fig:fig3} that the designed functional ODE observer works well for the given system.
\begin{figure}[h]
	\centering
	\begin{center}
\includegraphics[scale=0.60]{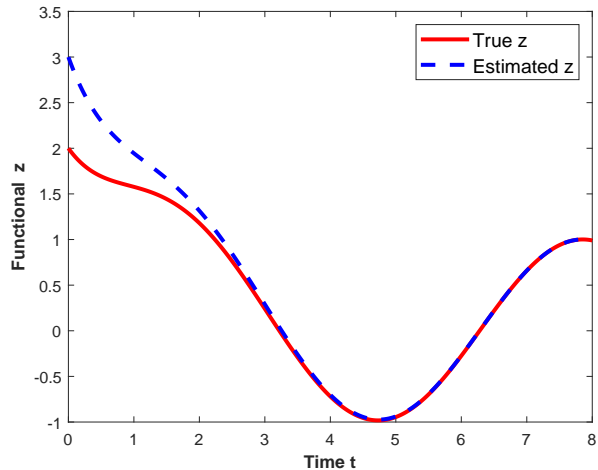}
	\caption{Time response true and estimated $z$.}
			\label{fig:fig3}
		\end{center}
	\end{figure}
	
\end{example}

\section{CONCLUSION}\label{sec:conc}
We extended the results of \cite{jaiswal2021existence} to obtain new findings
pertaining to the functional ODE observer design problem for
a comprehensive class of linear descriptor systems. Though it
has been shown that the derived conditions are
not necessary for the existence of functional ODE observers,
they are much milder than those obtained in the
earlier works.
Nevertheless, a topic of future research will be to close this gap and find simple conditions that are necessary and sufficient for the existence of functional ODE observers for linear descriptor systems.

\end{document}